\pdfoutput=1 
\documentclass[11pt]{amsart}
\usepackage{amsmath,amssymb}
\usepackage{graphicx,color}

\newtheorem{thm}{Theorem}[section]

\newtheorem{rem}[thm]{Remark}

\def\qed{\hfill $\Box$}
\def\proof{\noindent {\sl Proof} :\;  }

\newcommand{\B}{\mathcal{B}}
\newcommand{\A}{\mathcal{A}}
\newcommand{\K}{\mathcal{K}}

\newcommand{\R}{\mathbb{R}}

\def\qed{\hfill $\Box$}
\def\proof{\noindent {\sl Proof} :\;  }

\def\rd{\partial}

\begin{document}
\title[Bifurcation of parabolic map-germs]{Bifurcation of plane-to-plane map-germs 
with corank two of parabolic type }
\author[T.~Yoshida]{Toshiki Yoshida}
\address[T.~Yoshida]{Department of Mathematics, 
Graduate School of Science,  Hokkaido University,
Sapporo 060-0810, Japan}
\email{toshiki@mail.sci.hokudai.ac.jp}
\author[Y.~Kabata]{Yutaro Kabata}
\address[Y.~Kabata]{Department of Mathematics, 
Graduate School of Science,  Hokkaido University,
Sapporo 060-0810, Japan}
\email{kabata@math.sci.hokudai.ac.jp}
\author[T.~Ohmoto]{Toru Ohmoto}
\address[T.~Ohmoto]{Department of Mathematics, 
Graduate School of Science,  Hokkaido University,
Sapporo 060-0810, Japan}
\email{ohmoto@math.sci.hokudai.ac.jp}
%
%
\maketitle
\begin{abstract} 
There is a unique $\A$-moduli stratum of plane-to-plane germs 
which forms an open dense subset in 
the $\K$-orbit of $I_{2,3}: (x^2+y^3, xy)$.  
We describe explicitly the bifurcation diagram of 
its topologically $\A_e$-versal unfolding. 
Two geometric applications to parabolic crosscaps and parabolic umbilic 
are presented. 
\end{abstract}

\section{Introduction}

The bifurcation diagram of a family of smooth functions or mappings  
takes a fundamental role in {\it Catastrophe Theory} -- 
by definition it is the locus in the parameter space  at which 
the corresponding function fails to be  {\it structurally stable}; 
Along the locus, qualitative changes of the function occur.  
In this paper, we deal with 
generic families of plane-to-plane maps with at most $3$-parameters 
within the $\A$-classification theory 
($\A$ denotes the equivalence of map-germs via 
the action of diffeomerphism-germs of source and target). 

Bifurcation diagrams for $\A$-types of corank one germs can be found in several literatures;  
Arnold-Platonova \cite{Arnoldbooklet, Arnoldency}, 
Rieger \cite{Rieger2, Rieger3}, 
Gibson-Hobbs \cite{GH} and Aicardi-Ohmoto \cite{OA}. 
In contrast, for  corank two germs (Table \ref{list2} below), 
there had been very few known diagrams until quite recently.  
The case of deltoid is easy, while other cases may require 
a lot of computations; 
The bifurcation diagrams for types sharksfin 
and  odd-shaped sharksfin were first presented rigorously 
in  \cite{Yoshida, YKO}, 
that will partly be summarized in \S 3 (Fig.\ref{I22}, \ref{bifurcation_os}). 
The remaining case in codimension $3$ 
is the {\it $\A$-moduli stratum} in the $\K$-orbit of 
the germ $I_{2,3}: (x^2+y^3, xy)$; 
the $\K$-orbit has a unique topological $\A$-type 
(Rieger-Ruas \cite{RR},  Gaffney-Mond \cite{GM}). 
\begin{figure}
\begin{center}
\includegraphics[width=2.5cm, clip]{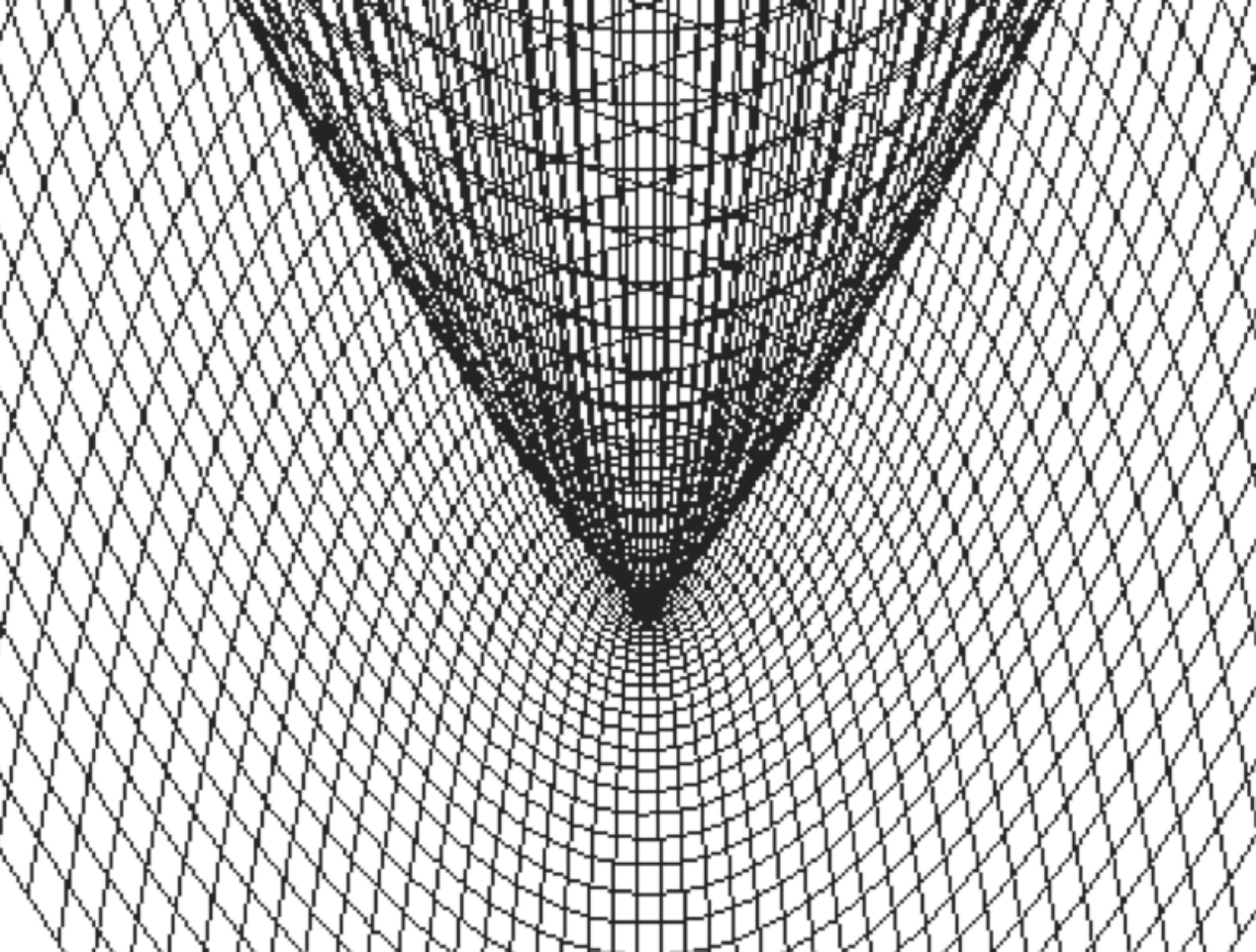} \quad \quad 
\includegraphics[width=3cm, clip]{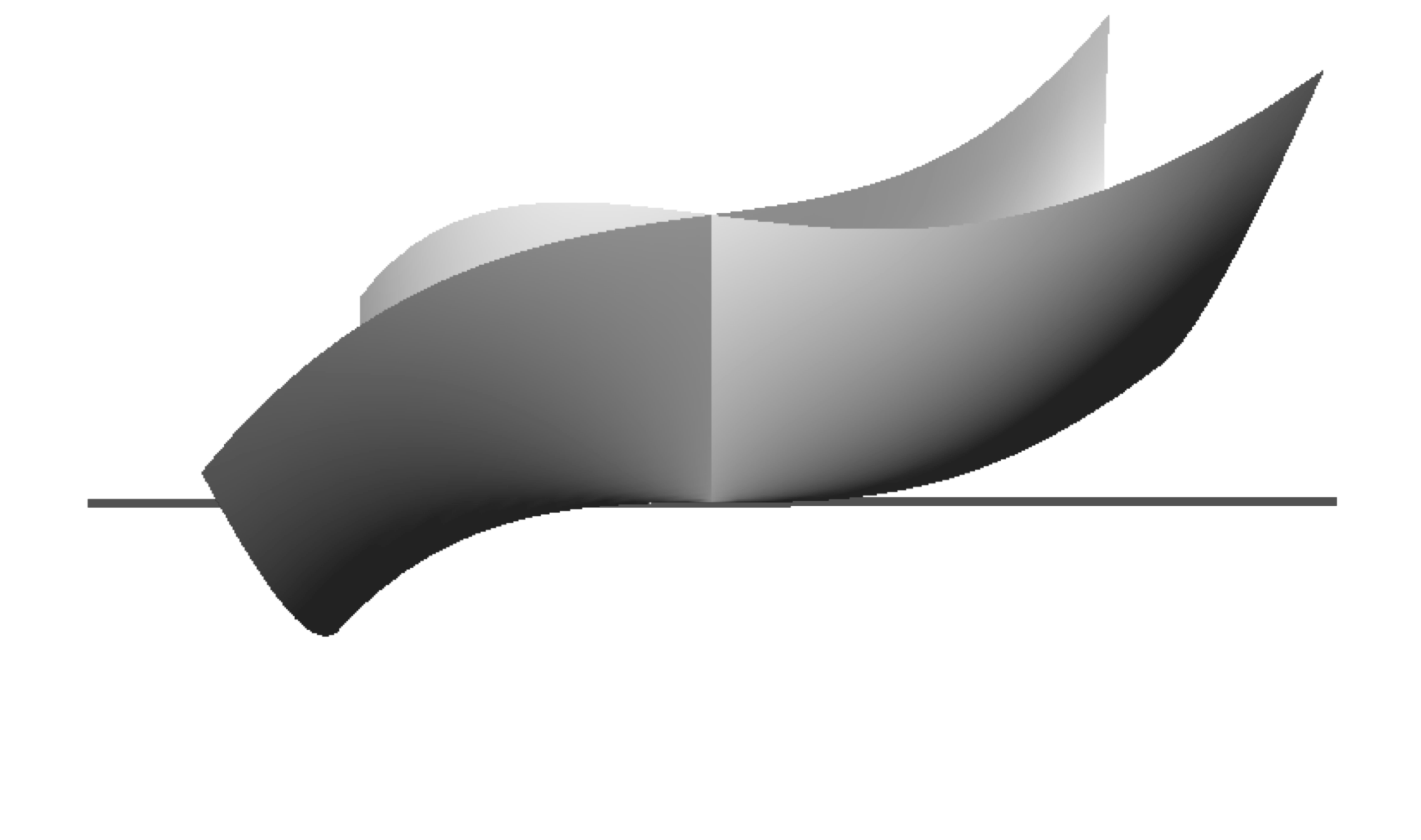} 
\end{center}
\caption{\small  Apparent contour of $I_{2,3}$ and projection of parabolic crosscap}
\label{proj}
\end{figure}
Our main purpose  is to compute and describe 
explicitly  its bifurcation diagram 
 in the same way as  \cite{Yoshida, YKO} (Theorem \ref{thm1} and Fig.\ref{I23}). 
Note that the structure of nearby $\K$-orbits 
has been known in 70's 
by Lander \cite[\S 5.5]{Lander}. 
In this paper 
we go further to find the finer structure of 
local and multi-$\A$-types  
appearing in the topologically $\A_e$-versal unfolding of type $I_{2,3}$.

\begin{table}
\begin{center}
$$
\begin{array}{c| l | l }
\mbox{cod}\;   &\;\; \mbox{type}  &\;\; \mbox{normal form}\\
\hline
\hline
2&I_{2,2}^{1,1}\, (\mbox{sharksfin}) & (x^2+y^3, y^2+x^3) \\
&II_{2,2}^{1,1}\, (\mbox{deltoid}) & (x^2-y^2+x^3, xy) \\
\hline
3&I_{2,2}^{2,1} \, (\mbox{odd-shaped sharksfin}) & (x^2+y^5, y^2+x^3) \\
&I_{2,3}\, (\mbox{$\A$-moduli}) & (x^2+y^3+\alpha xy^2 +\beta y^4+ \cdots, xy) \\
\end{array}
$$
\end{center}
\caption{\small Corank two germs with $\A_e$-codim. $\le 3$. }
\label{list2}
\end{table}

In the last half of this paper, we discuss applications to 
`parabolic objects' in several geometric settings. 
The key idea of our approach is as follows: 
We first reformulate the problem 
in terms of $\A$-singularity types arising 
in some family of plane-to-plane maps naturally associated to the setting, 
and then  embed the family 
into our $C^0$-versal unfolding of $I_{2,3}$. 
That may yield a $2$-dimensional section of our bifurcation diagram 
in the parameter $3$-space, 
from which 
one can deduce some topological nature of 
bifurcations of the parabolic object 
under consideration. 
For instance, the germ of type $I_{2,3}$ 
naturally appears in the parallel projection of 
a {\it parabolic crosscap} to the plane along 
the tangent line at the crosscap point (Fig.\ref{proj}, right). 
Our result provides a new insight into  
differential geometry of parabolic crosscaps in $\R^3$
(cf. Nu\~no-Ballesteros and Tari \cite{NT},  Oliver \cite{Oliver}). 
In geometric optics or symplectic geometry, 
the singularity of $I_{2,3}$ is realized as a {\it planar caustics of parabolic umbilic type $D_5$}, 
which is one of most favorite singularities in 
R. Thom  \cite{Thom}. 
In Appendix, we apply our result to this setting 
and discuss 
generic $2$-parameter `perestroikas' of planar caustics. 

The authors thank the organizers of Japanese-Brazilian workshop 
on singularities (RIMS, 2013); it was a nice opportunity to discuss the problem 
dealt in this paper with several experts. 
The third author is partly  supported  by JSPS grant no.23654028.

\section{Recognition of $\A$-types of corank one germs}
\subsection{$\A$-classification}
$C^\infty$ map-germs $f, g: \R^2, 0 \to \R^2,0$ are  {\it $\A$-equivalent} 
if there is a pair $(\sigma, \tau)$ of diffeomorphism germs 
of source and target planes 
at the origins so that $g=\tau \circ f \circ \sigma^{-1}$. 
We say $f$ is $\A$-simple if 
the cardinality of nearby $\A$-orbits is finite 
(otherwise, $f$ belongs to an $\A$-modulus, i.e., 
some continuous family of $\A$-orbits). 
The {\it corank} of $f$ means $\dim \ker df(0)$. 
Let $T\A_e.f$ denote the extended $\A$-tangent space of $f$, then 
the {\it $\A_e$-codimension of $f$} is defined by $\dim_\R \theta(f)/T\A_e.f$; 
that is the smallest number of parameters required 
for constructing its $\A_e$-versal unfolding 
(such an unfolding is called to be ($\A_e$-)miniversal). 
Equivalently, a singularity type of $\A_e$-codimension $\le r$ 
means that it generically appears in $r$-parameter families of maps. 
In particular, 
$f$ is a {\it stable germ} if and only if its $\A_e$-codimension is $0$. 
As for the $\A$-classification of plane-to-plane germs, see 
Rieger \cite{Rieger}  (for corank $1$ germs) and Rieger-Ruas \cite{RR} (for corank $2$ germs). 

\begin{figure}
\begin{center}
\includegraphics[width=13cm, clip]{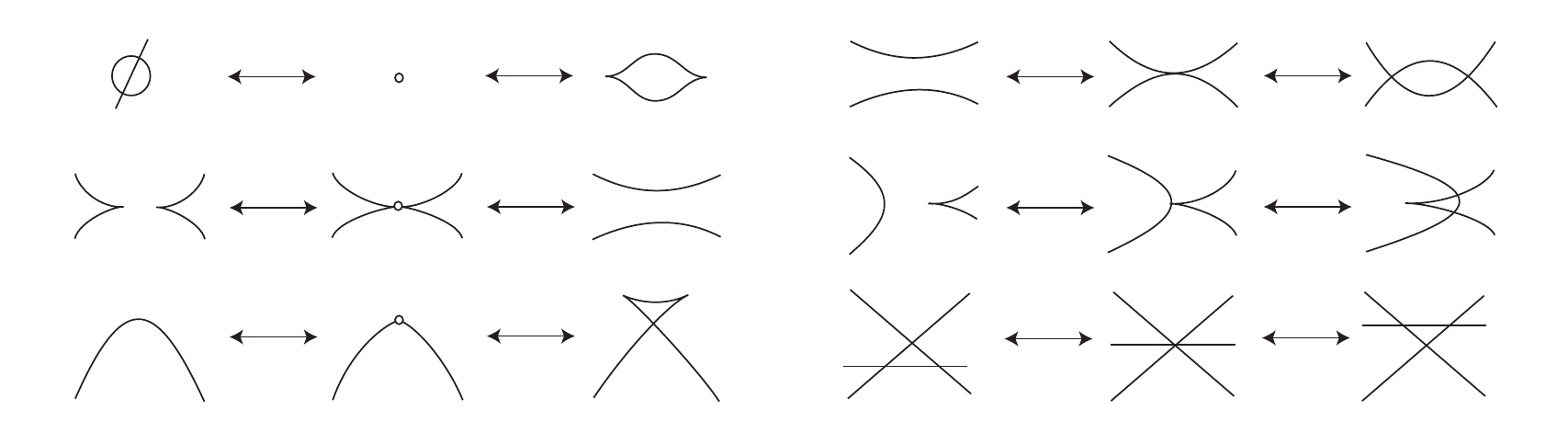}
\end{center}
\caption{\small  Apparent contour of germs with $\A_e$-codimension $1$ \cite{OA}: 
local singularity types are of  
lips, beaks, swallowtail (left), and  multi-singularity types are of tacnode folds, cusp+fold, triple folds (right). }
\label{codim1}
\end{figure}

\begin{figure}
\begin{center}
\includegraphics[width=13cm, clip]{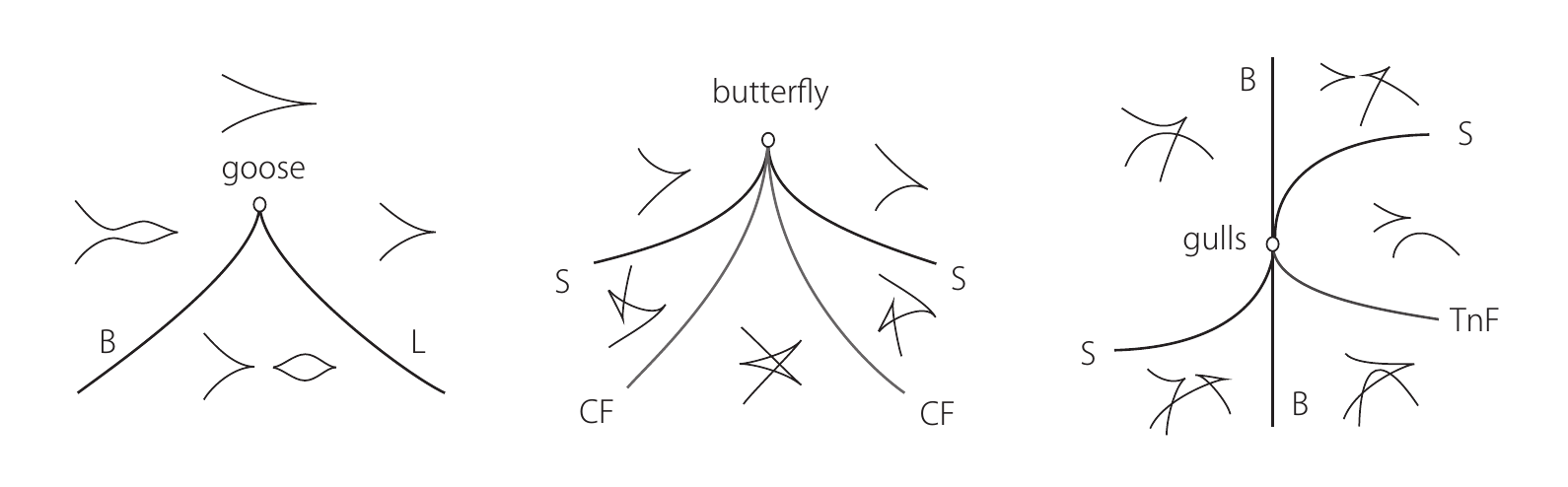}
\end{center}
\caption{\small  Bifurcation diagrams of  singularities with corank one of 
$A_e$-codimension $2$:  
goose $(x, y^3+ x^3y)$, butterfly $(x,xy+y^5+y^7)$, gulls $(x,xy^2+y^4+y^5)$.  
}
\label{codim2}
\end{figure}

Let $f:\R^2,0 \to \R^2,0$ be a germ of $\A_e$-codimension $s$, 
and $F$ an $\A_e$-miniversal unfolding of $f$. 
Take a good representative 
$F: U \times W \to \R^2 \times W$, where $U \subset \R^2$ are $W \subset \R^s$ 
are sufficiently small open neighborhoods of origins, 
and consider a $C^\infty$ map $F_\lambda: U \to \R^2$ for each $\lambda \in W$. 
For general $\lambda$, 
the map $F_\lambda$ is stable, that is, 
$F_\lambda$ has only singularities of type {\it fold},  {\it cusp} and 
{\it double folds} (bi-germ). 
The {\it bifurcation diagram} $\B_F$ is defined to be 
the locus consisting of $\lambda$ so that 
$F_\lambda$ has {\it unstable} (mono/multi)-singularities 
at some points in $U$. 
In particular, $\B_F$ is stratified 
according to local and multi-singularity types of germs 
with $\A_e$-codimension less than $s$. 
For instance, all local and multi singularity types 
of $\A_e$-codimension $1$ are presented in Fig.\ref{codim1}, 
and  local singularity types of codimension $2$ 
are depicted in Fig.\ref{codim2} (corank one) and Fig.\ref{I22} (corank two); 
Besides, there are 15 types of multi-singularities of codimension $2$ 
(some combinations of codimension one singularities), see \cite{OA}. 

\subsection{$\A$-recognition and geometric criteria} 
To find an explicit equation for each stratum in $\B_F$, 
there is a useful tool \cite{Saji, Kabata}. 
Given a map-germ $f=(f_1, f_2): \R^2, 0 \to \R^2,0$ of corank one,  
we want to determine which $\A$-type the germ $f$ belongs to. 
Let us take 
\begin{itemize}
\item the  Jacobian 
$\lambda(x,y):=\frac{\partial (f_1, f_2)}{\partial (x, y)}$
\item  arbitrary non-zero vector field $\eta=\eta_1(x,y)\frac{\rd}{\rd x}+\eta_2(x,y)\frac{\rd}{\rd y}$ 
near the origin which spans $\ker df$ 
at any singular points $\lambda=0$: 
\end{itemize}
We put $\eta^kg=\eta(\eta^{k-1}(g))$ for any function $g(x,y)$. 
Additionally, if the Hessian matrix $\mbox{H}_\lambda$ of $\lambda$ at $0$ has rank one, 
let $\theta$ be a vector field so that $\theta(0)$ spans $\ker \mbox{H}_\lambda (0)$. 
Then a geometric characterization of each $\A$-type 
with $\A_e$-codimension $\le 2$ 
is described  
in terms of $\lambda$ and $\eta$ (and $\theta$) as in Table \ref{criteria}
(see \cite{Kabata}, for $\A$-types with higher codimension). 

\begin{table}
\begin{center}
\begin{tabular}{ll}
\hline 
fold& 
 $\eta \lambda (0) \neq 0$\\
cusp & 
$ d \lambda (0) \neq 0$, $\eta \lambda (0) = 0$, $\eta^2 \lambda (0) \neq 0$\\
\hline
swallowtail &
$d\lambda (0) \neq 0$, $ \eta \lambda (0) = \eta^2 \lambda (0) =0$, 
$ \eta^3 \lambda (0) \neq 0$\\
lips  &
$d \lambda (0) = 0$, $ \det \mbox{H}_\lambda (0) > 0$
\\
beaks &
$d \lambda (0) = 0$, $ \det \mbox{H}_\lambda (0) <0 $, ${\eta} ^2 \lambda (0) \neq 0$
\\
\hline 
butterfly  & 
$d \lambda (0) \neq 0$, $\eta \lambda (0) = \eta^2 \lambda (0) =\eta^3 \lambda (0) =0$, 
$ \eta^4 \lambda (0) \neq 0$
\\
gulls &
$d \lambda (0) = 0$, $\det \mbox{H}_\lambda (0)< 0$, 
${\eta} ^2 \lambda (0) =0,  {\eta} ^3 \lambda (0) \neq 0$
\\
goose   & 
$d \lambda (0) = 0$, $\mbox{rk}\,  \mbox{H}_\lambda (0) =1$, 
$ {\eta} ^2 \lambda (0) \neq 0, \theta^3 \lambda (0) \neq 0.$\\
\hline \\
 \end{tabular}
 \end{center}
 \caption{\small Criteria (of jets) of $\A$-orbits}
 \label{criteria}
\end{table}

\section{Corank two singularities}

\subsection{Sharksfin and odd-shaped sharksfin}
\begin{figure}
\begin{center}
\includegraphics[width=12.5cm, clip]{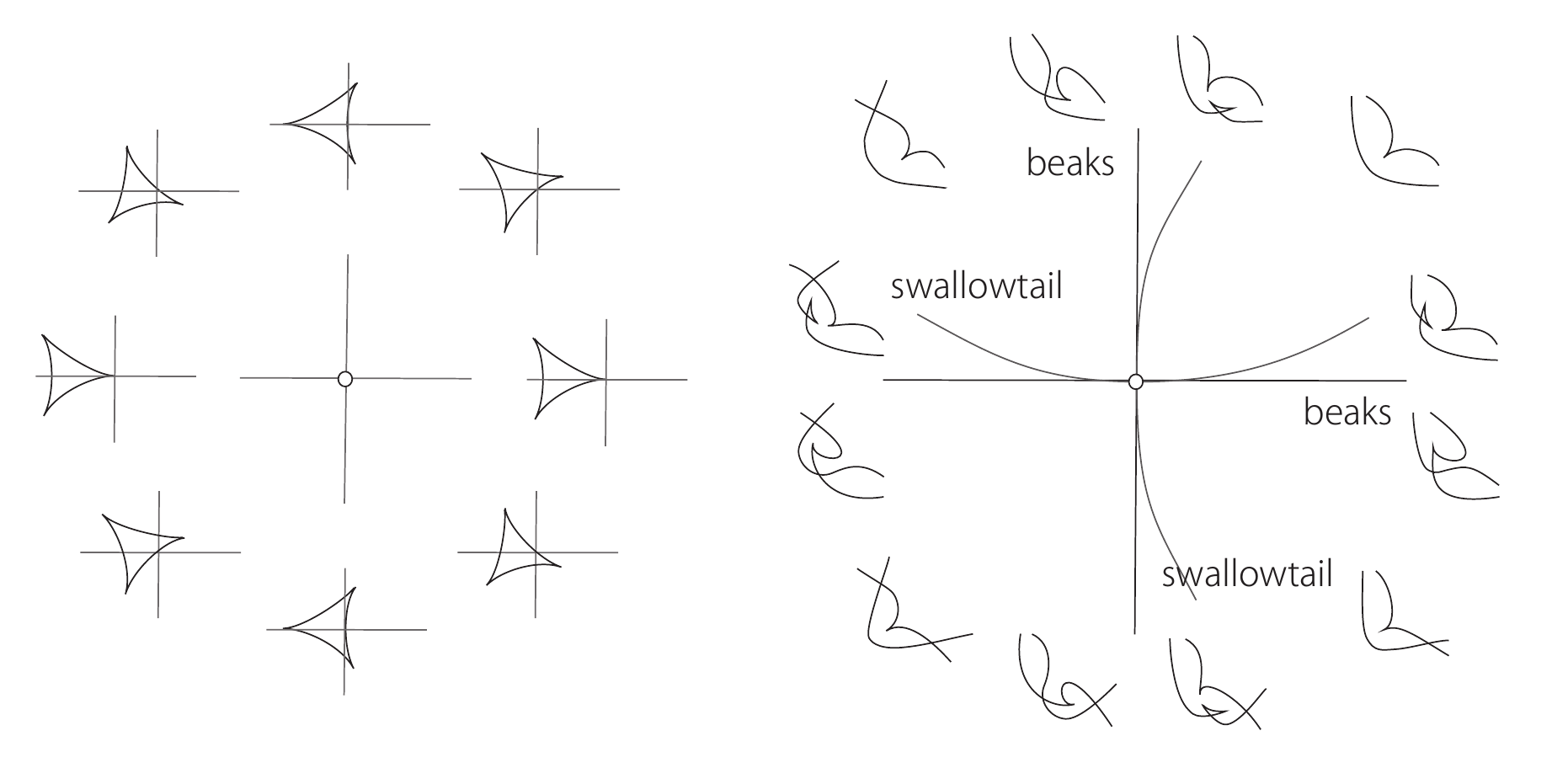} 
\end{center}
\caption{\small Bifurcation diagrams of singularities with corank two of 
$A_e$-codimension $2$: deltoid 
$(x^2-y^2+x^3, xy)$ and sharksfin  $(x^2+y^3, y^2+x^3)$ \cite{GH, YKO}}
\label{I22}
\end{figure}

 In Rieger-Ruas \cite{RR}, $\A$-simple germs of 
corank $2$ have been classified. 
As seen before (Table \ref{list2}), 
there are four types in codimension $\le 3$, 
one of which is not $\A$-simple, that is the moduli of type $I_{2,3}$.  
The case of deltoid is easy: 
the bifurcation diagram consists only of the origin in the parameter plane, i.e., 
the singularity type has only adjacencies of 
fold and cusps, and no other local and multi-singularities; 
any small perturbation of this type produces a `deltoid-shaped' apparent contour 
with three cusps (Fig.\ref{I22}, left). 
On the other hand, 
the bifurcation diagrams of sharksfin and odd-shaped sharksfin 
had been unclear for a long time \cite{Hawes, Yoshida, YKO}.

\begin{figure}
\begin{center}
\includegraphics[width=8.5cm, clip]{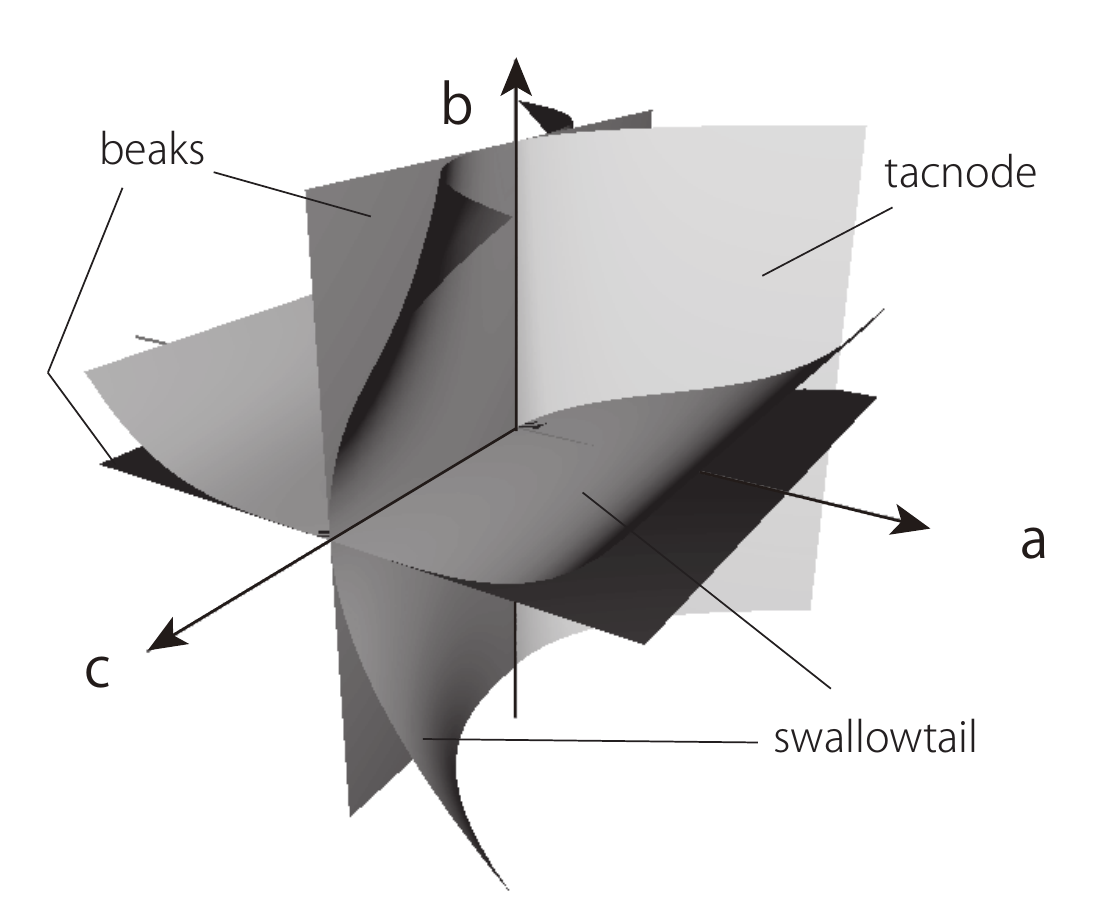}
\end{center}
\caption{\small Bifurcation diagrams of  
 odd-sharksfin \cite{YKO}}
\label{bifurcation_os}
\end{figure}

\subsubsection{\bf Sharksfin} 
Let us consider the following miniversal unfolding of  $I_{2,2}^{1,1}$: 
\begin{equation}\label{versal_sharksfin}
F(x,y,a,b)=(x^2+y^3+ay, y^2+x^3+bx).
\end{equation}
In the parameter $ab$-plane,  the bifurcation diagram $\B_{F}$ consists of 
four smooth curve-germs at the origin 
(for the swallowtail, it is computed up to degree $9$ in \cite{Yoshida}): 

\

\begin{center}
\begin{tabular}{ll}
{\bf  Beaks}: &
$a=0$ {\rm and} $b=0$ \\
{\bf Swallowtail}: &
$a=\frac{1}{16}b^4+\frac{3}{32}b^9+o(9)$ {\rm and} 
$b=\frac{1}{16}a^4+\frac{3}{32}a^9+o(9)$. 
\end{tabular}
\end{center}

\subsubsection{\bf Odd-shaped sharksfin} 
Let us consider the following miniversal unfolding  of   $I_{2,2}^{2,1}$: 
\begin{equation}\label{versal_oss}
F(x,y,a,b,c)=(x^2+y^5+cy^3+ay, y^2+x^3+bx).
\end{equation}
The 3D picture of $\B_{F}$ is drawn in Fig.\ref{bifurcation_os}. 
Although it is too hard to find 
an explicit form of the defining equation for the swallowtail stratum, 
we know how to draw it from the information about 
the bifurcation diagrams of types sharksfin and gulls.

\

\begin{center}
\begin{tabular}{ll}
{\bf  Beaks}: & $a=0$ {\rm and} $b=0$ \\
{\bf  Swallowtail}: & {\rm two smooth surfaces tangent to $ab=0$ along $c$-axis} \\
&{\rm with $4$-point contact, one of which is tangent to $a=0$}\\
&{\rm along the $b$-axis with $3$-point contact. }\\
{\bf Tacnode}: & $4a=c^2\; (c<0)$\\
{\bf Gulls}: & {\rm $b$-axis}
\end{tabular}
\end{center}

\subsection{Bifurcation diagram of type $I_{2,3}$}
It is known  that the germ\footnote{
Usually Mather's  notation $I_{2,3}$ is used for the $\K$-orbit, 
but in this paper we use it for this particular map-germ.}  
$$I_{2,3}:(x,y) \mapsto (x^2+y^3, xy)$$ 
is $\A$-finite, and is not $\A$-simple -- 
it belongs to a moduli stratum 
of $\A$-orbits of the form 
$(x^2+y^3+\alpha xy^2+\beta y^4+\cdots, xy)$ 
with the modality $\ge 2$, see Rieger-Ruas \cite{RR}. 
Furthermore, 
as noted in Gaffney-Mond \cite[Ex.5.11]{GM}, 
any $\A$-finite germ contained in the $\K$-orbit is obtained by 
adding some higher terms to this germ, and hence 
the germ is topologically $\A$-equivalent to $I_{2,3}$ 
by a theorem of J. Damon, 
i.e., the $\K$-orbit has a unique topologically $\A$-type 
in its open dense subset. 
In other words, 
in the larger parameter space of 
an $\A_e$-miniversal unfolding of the $\A$-finite germ $I_{2,3}$, 
the bifurcation diagram is Whitney regular 
along the $\A$-moduli stratum. 
Therefore for our purpose,  
it suffices to take the following form of 
an unfolding of $\A$-type $I_{2,3}$ 
which corresponds to a normal slice to the stratum: 
\begin{equation}\label{versal_I23}
G(x,y,a,b,c)=(x^2+y^3+ax+by+cy^2, xy).
\end{equation}
Namely, even if we add some terms $xy^2, y^4, \cdots$ to the first component, 
the unfolding remains to be transverse to the stratum in the space of all germs (or jets), 
and thus the bifurcation diagram is topologically the same 
by Thom's isotopy lemma. 

Note that $G$ is just a stable unfolding of $\K$-type of $I_{2,3}$: 
The structure of nearby $\K$-orbits was well investigated 
in  Lander's paper \cite{Lander}, 
in which  the loci of swallowtail and butterfly of $\B_G$ 
are presented in an explicit form. 
Our first aim  is to describe precisely 
all the strata of the bifurcation diagram $\B_G$ 
for local and multi $\A$-types. 

\begin{thm}\label{thm1}
Let $G$ be the topologically $\A_e$-versal unfolding (\ref{versal_I23}) of $I_{2,3}$. 
Then the bifurcation diagram $\B_G$ consists of 
three components  
corresponding to types beaks-lips, swallowtail and cusp+fold, 
as drawn in Fig.\ref{I23}. 
Each stratum is explicitly parametrized  as follows. 
For the complement to $\B_G$, 
apparent contours of corresponding stable maps 
are drawn in Fig.\ref{section_I23+}. 

\

\begin{center}
\begin{tabular}{ll}
{\bf Sharksfin} &  $c$-axis with $c>0$\\
{\bf Deltoid} &  $c$-axis with $c<0$\\
{\bf  Beaks/Lips}: &   
the surface with $A_3$-singularity at the origin whose \\
&  double point curve  is the locus of sharksfin; parametrized by \\
&  $ (a,b,c)=(\pm 4y\sqrt{c+3y}, -y(4c+9y), c)$ 
\vrule width 0pt height 12pt depth 10 pt\\
{\bf Goose}: 
&  the cuspidal edge of the beaks/lips locus parametrized by\\
&  $(\pm \frac{8}{9\sqrt{3}}c^{3/2}, \frac{4}{9}c^2, c)$ with $c>0$ 
\vrule width 0pt height 12pt depth 10 pt\\
{\bf  Swallowtail}: &   the surface which contains 
the $a$-axis and the sharksfin locus; \\
& parametrized by  \\
& $ (a,b,c)=(\pm y(4c+15y)(c+4y)^{-1/2}, -2y(2c+5y), c)$
\vrule width 0pt height 12pt depth 10 pt \\
{\bf Butterfly}: & the cuspidal edge of the swallowtail locus parametrized by \\
& $(\pm \frac{1}{\sqrt{5}}c^{3/2}, \frac{2}{5}c^2, c)$ with $c>0$
\vrule width 0pt height 12pt depth 10 pt\\
{\bf Cusp+Fold}: 
& the surface with $A_3$-singularity at the origin whose\\
&  cuspidal edge  is the butterfly locus;  parametrized by \\
& $ (a,b,c)=(\pm \sqrt{-y}(3c+5y), \frac{1}{4}(c^2-6cy-15y^2), c)\;\; (y<0)$
\vrule width 0pt height 12pt depth 10 pt\\
\end{tabular}
\end{center}
\end{thm}

\begin{rem}{\rm 
It is remarkable that 
the odd-shaped sharksfin has the adjacency of gulls, 
but not butterfly and goose, 
while 
the type of $I_{2,3}$ has the adjacencies of butterfly and goose, but not gulls. 
That completes the adjacency diagrams 
among $\A$-orbits of $\A_e$-codimension $\le 3$ as in Fig.\ref{adjacency} 
(it corrects \S 3 of \cite{RR}). We follow the notation used in \cite{Rieger} for each $\A$-type. 

\begin{figure}[htb]
\begin{center}
\includegraphics[width=13cm, clip]{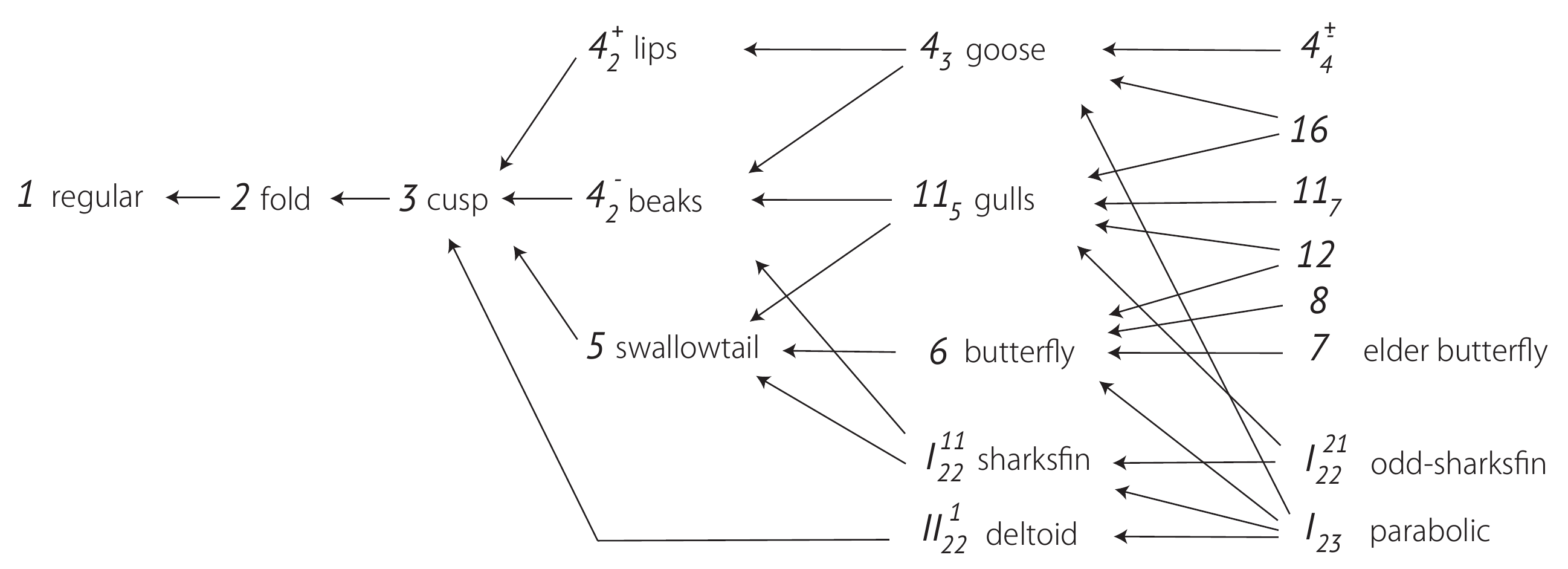}
\end{center}
\caption{\small Adjacencies of local singularities up to $\A$-cod. $\le 5$. }
\label{adjacency}
\end{figure}
}
\end{rem}

\begin{figure}
\begin{center}
\includegraphics[width=9cm, clip]{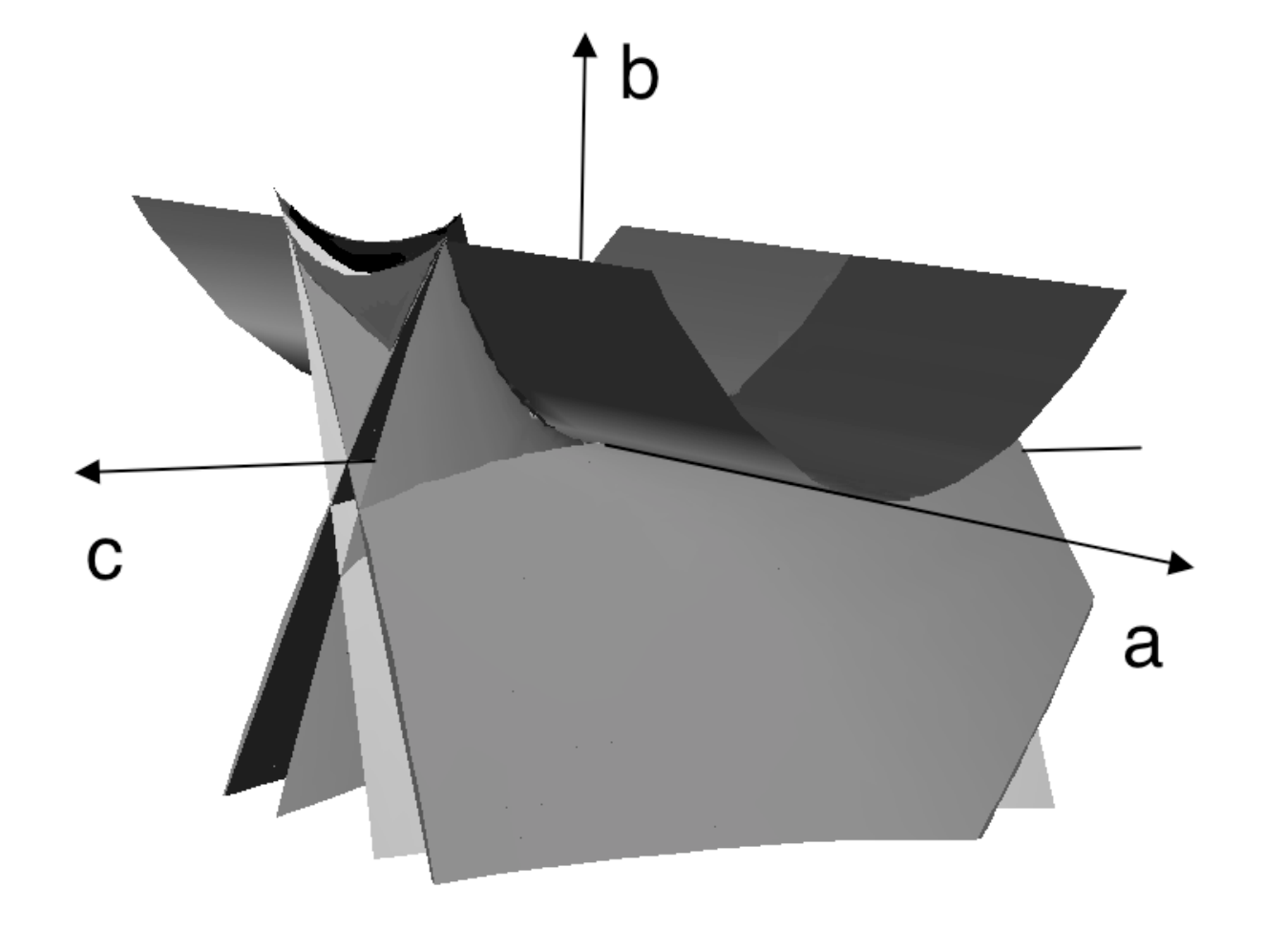}\\
\includegraphics[width=12.5cm, clip]{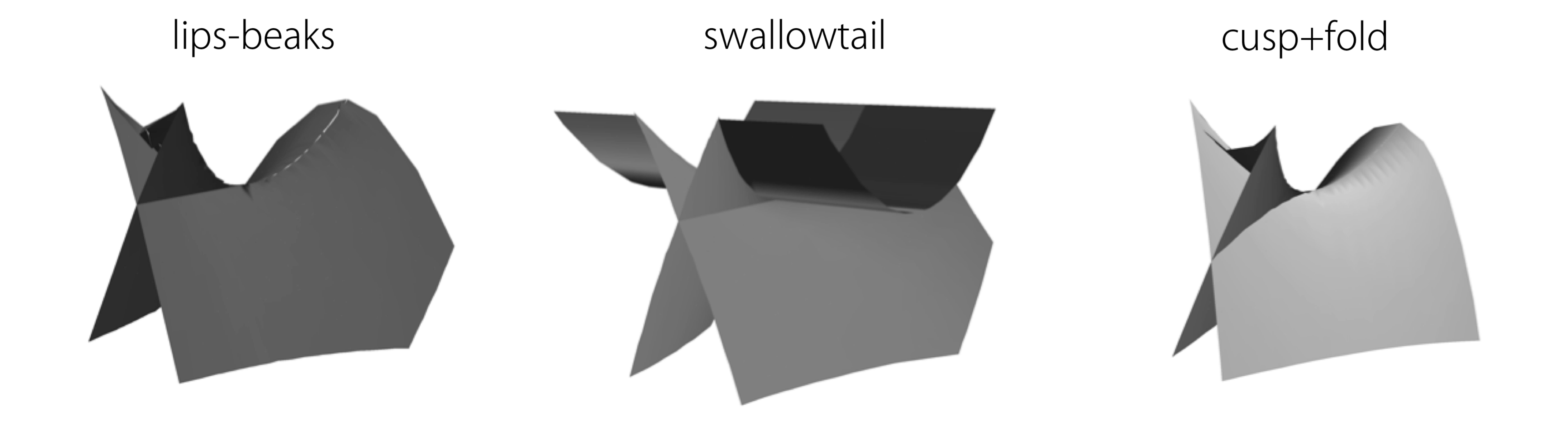}
\end{center}
\caption{\small 
The bifurcation diagram in the parameter $3$-space consists of three `doggies' components, 
beaks-lips locus (left), swallowtail locus (center), cusp+fold locus (right), 
and a `doggie's tail' component (one half of the $c$-axis with $c<0$) 
which is the deltoid locus. 
The picture of the swallowtail locus coincides with Fig.3 in Lander \cite{Lander} 
(but from a different viewpoint). 
 }
\label{I23}
\end{figure}

\begin{figure}
\begin{center}
\includegraphics[width=13cm, clip]{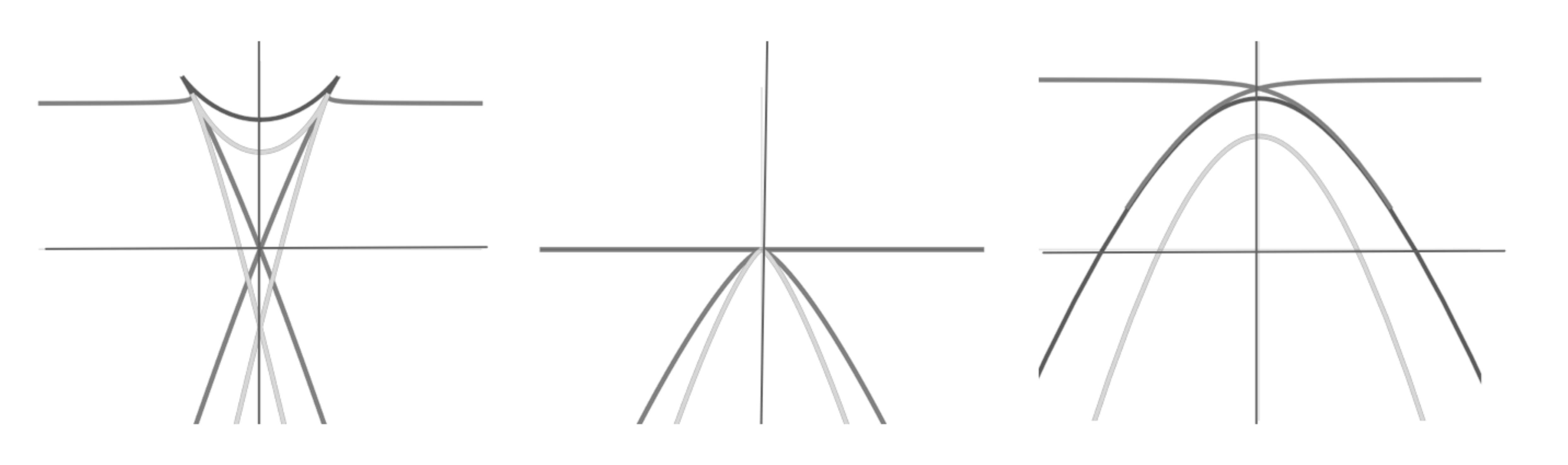}
\end{center}
\caption{\small Sections of the bifurcation diagram 
with the plane $c=const.$ (positive, zero, negative, from left to right). 
The loci of 
beaks and swallowtail are very close to each other in a large part, 
where they are depicted to be overlapped. 
See Fig. \ref{section_I23+}. }
\label{section_I23}
\end{figure}

\begin{figure}
\begin{center}
\includegraphics[width=13cm, clip]{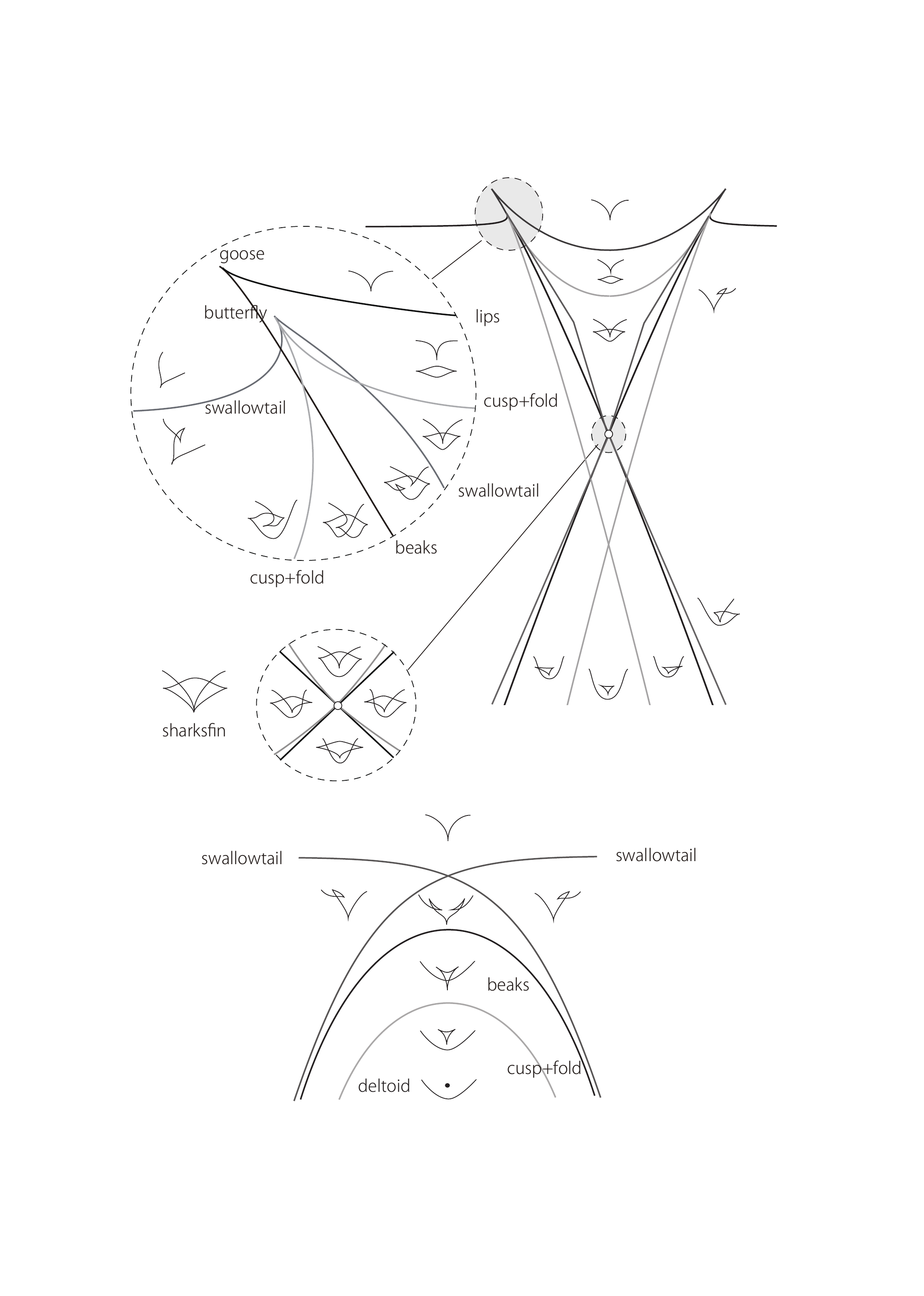}
\end{center}
\caption{\small   
Bifurcation of 
apparent contours of map-germs corresponding to sections in Fig. \ref{section_I23} (left and right).  }
\label{section_I23+}
\end{figure}

\newpage

\proof 
To obtain above parameterizations 
we use the criteria in Table \ref{criteria}. 
For the unfolding given by (\ref{versal_oss}), 
put 
$$\textstyle \lambda=
\left|
\begin{array}{cc}
\;2x+a &\;\; \;\;3y^2+2cy+b\;\\
\;y&\;\;\;\;x\;
\end{array}
\right|, 
$$
and take 
$$\eta=-x \frac{\rd}{\rd x}+y\frac{\rd}{\rd y}$$
where $(x,y)\not=(0,0)$, 
otherwise $\eta=-(3y^2+2cy+b)\frac{\rd}{\rd x}+(2x+a)\frac{\rd}{\rd y}$. 

To see the loci of sharksfin and deltoid is easy: $x=y=a=b=0\; (c\not=0)$. 
The sharksfin (resp. the deltoid) corresponds to $c>0$ (resp. $c<0$). 

The beaks/lips locus is defined by 
$$\textstyle 
\lambda=0, \; \frac{\rd}{\rd x}\lambda=a+4x=0, \; 
\frac{\rd}{\rd y}\lambda=-b-4cy-9y^2=0.$$
Solve $a, b$ in $x, y, c$ and substitute them into $\lambda=0$, 
then we have  $x=\pm y\sqrt{c+3y}$ 
(thus in $xyc$-space it forms a crosscap $x^2=y^2(c+3y)$). 
This equality gives the above parametrization of $a=a(c,y)$ and $b=b(c,y)$. 
The picture of this locus is depicted in Fig.\ref{I23} (left): 
It has a transverse self-intersection along the sharksfin locus. 
The defining equation  is given by 
$$243 a^4+ \left(256 c^3-864 b c\right)a^2+768 b^3-256 b^2 c^2=0.$$
The strata of type lips and the beaks are separated by the goose curve: 
From the criteria (Table \ref{criteria}), 
we check $\mbox{H}_\lambda=0$ additionally; 
that implies $y=-\frac{2}{9}c\; (c>0)$,  
and hence the above parametrization is obtained. 
The goose curve is actually  the cuspidal edge of the singular surface, 
which looks like the first doggies ``ears". The lips part is  the ``head".  
The gulls type does not appear, 
since $\eta^2\lambda=0$ implies $x=y=a=b=0$. 

The swallowtail locus is defined by 
$$\lambda=\eta\lambda=\eta^2\lambda=0.$$
It is solved as follows. 
The computation is essentially the same as in Lander \cite[\S 5.5]{Lander}, since $G$ is a stable unfolding of $\K$-type of $I_{2,3}$. 
By $\lambda=\eta\lambda=0$, we have
$$\textstyle a=\frac{1}{x}(-2x^2+by+2cy^2+3y^3), \;\; 
b=-\frac{1}{y}(x^2+3cy^2+6y^3),$$ 
unless $xy=0$, 
and then  $\eta^2\lambda=0$ leads to $x=\pm y\sqrt{c+4y}$ 
(thus in $xyc$-space it forms a crosscap). 
This yields the above parameterization. 
Notice that the limit as $c, y \to 0$ is the $a$-axis; 
In fact, if $x=0$ or $y=0$, three equations imply $x=y=b=c=0$, 
hence the $a$-axis belongs to this locus. 
The defining equation of the locus is 
$$80 a^4 \left(8 b-3 c^2\right)-8 a^2 \left(285 b^2 c-192 b c^3+32 c^5\right)+b^2\left(45 b-16 c^2\right)^2=0, $$
and the picture is Fig.\ref{I23} (center). 
It has a transverse self-intersection along the sharksfin locus: 
The beaks and the swallowtail loci have $4$-point contact 
along this half line, 
that is verified by the bifurcation diagram of the sharksfin. 
For the butterfly, we add one more equation 
$\eta^3\lambda=-24y^3(c+5y)=0$, that implies  $y=-c/5\; (c>0)$; 
The butterfly curve also looks ``ears" of the second doggie, i.e., 
it is the cuspidal edge of the locus of swallowtail. 

The locus of cusp + fold (bi-germ) must appear,  since it is adjacent to the butterfly. 
The equations are: 
\begin{eqnarray*}
&&xy=XY, \; ax+x^2+by+cy^2+y^3=aX+X^2+bY+cY^2+Y^3, \\
&&\lambda(x,y)=\lambda(X,Y)=0, \; \eta\lambda(x,y)=0. 
\end{eqnarray*}
We can eliminate $a, b$ by the third and fourth equations 
$\lambda=0$; 
then the second equation  
gives $X=\pm y \sqrt{c+2(y+Y)}$, 
and hence the last equation eliminates $x$, then 
$Y=-\frac{1}{2}(c+3y)$. Thus we can express $a, b, x$ using variables $c, y$. 
The parametrization leads to  Fig.\ref{I23} (left).  
This third doggie also has ``ears" along the butterfly curve, as same as the second one does. 
The defining equation is 
\begin{eqnarray*}
&&
18225 a^8+14580 a^6 c^3-54 a^4 \left(1275 b^2 c^2-49 c^6\right)\\
  &&+108 a^2 \left(500 b^4   c-15 b^2 c^5-c^9\right) -\left(16 b^2-c^4\right) \left(c^4-25 b^2\right)^2=0. 
\end{eqnarray*}

In a similar way as seen above, 
direct computations show that 
no other multi-singularity strata appear. \qed

\section{Parallel projection of parabolic crosscaps}
\subsection{Projection of smooth surface in $3$-space} 
We begin with projecting a smooth surface to the plane. 
Let $M$ be a fixed smooth surface in $\R^3$, 
and $\Pi_\ell: \R^3 \to \R^2$ a linear projection with 
the kernel line $\ell \subset \R^3$, 
then the restriction  $\Pi_\ell|_M: M \to \R^2$ is called 
a {\it parallel projection} of $M$ along the direction $\ell$. 
When $\ell$ varies, 
it defines locally  a family of plane-to-plane maps 
with two parameters; 
V. I. Arnold \cite{Arnoldbooklet,  Arnoldency} and W. Bruce \cite{Bruce} classified 
singularities of the parallel projection for an appropriately generic surface 
up to the {\it $\A$-equivalence} defined by 
 local diffeomorphisms of $M$ and the target $\R^2$ (the screen). 
Such a generic surface is stratified according to 
the local singularity types of the projection. 
In particular, there are two major characteristic curves on $M$: 
\begin{itemize}
\item[-] 
The {\it parabolic curve} consists of points 
where the Gaussian curvature vanishes; 
At each point there is only one asymptotic line, 
and the parallel projection along the asymptotic line admits  
the lips and the beaks or more degenerate singularity.  
\item[-] 
The {\it flecnodal curve} consists of points where  
an asymptotic line has at least $4$-point contact with the surface; 
or equivalently, points where the Pick invariant vanishes; 
the parallel projection along the asymptotic line 
has the swallowtail or more degenerate singularity. 
\end{itemize}
Note that these two curves meet tangentially at 
some isolated points, called the {\it godrons}, 
where the projection admits the gulls singularity. 

Furthermore, not fixing a generic surface, 
we may consider 
a generic $1$-parameter family of embeddings 
$\iota_t: M \hookrightarrow \R^3$ ($t \in I$, where $I$ is an open interval); 
In relation with  an application to Computer Vision,  
J. Rieger \cite{Rieger3} studied  
singularities arising in the family $\Pi_\ell \circ \iota_t$ 
with three parameters $\ell$ and $t$, and showed that 
all singularity types of $\A_e$-codimension $\le 3$ 
arise generically. 

\begin{rem}{\rm 
In relation with projective geometry of surfaces, 
singularities in the {\it central projection from arbitrary viewpoint }
has been studied  by Arnold et al \cite{Arnoldbooklet,Arnoldency,Platonova, Goryunov}. 
For central projections of a moving surface, see Kabata \cite{Kabata}.  
The classification of singularities for central projections 
becomes slightly different from that for parallel projections. 
For instance, the goose singularity appears in parallel projection 
at some isolated parabolic points in a generic surface, 
while the singularity type always arises in central projection 
at any parabolic point, when viewing it 
from some special viewpoint lying on the asymptotic line 
(Note that the parallel projection corresponds to 
the central projection with the viewpoint at infinity). 
In this paper we only consider the parallel projection. 
}
\end{rem}

\subsection{Projection of crosscaps in $3$-space} 
As a generalization, 
J. West \cite{West} considered 
singularities of parallel projection 
$\Pi_\ell \circ \iota: M \to \R^2$ 
where $\iota: M \to \R^3$ is a smooth map having {\it crosscaps}. 
That is the (unique) locally stable singularity type, 
and if one take suitable local coordinates of source and target, 
the map-germ is written by $(y, xy, x^2)$. 
Since we are discussing 
the affine or flat geometry of the singular surface in $\R^3$, 
the ambient coordinate changes should be only affine transformation of $\R^3$, 
then the affine normal form is given by 
\begin{equation} \label{crosscap}
\iota(x,y)= (y, xy+g(x), x^2+\alpha y^2+\phi(x,y))
\end{equation}
for a constant $\alpha$, 
$g(x)=d_4x^4+h.o.t.$ and 
$\phi(x,y)=c_{03}y^3+c_{12}xy^2+\cdots $  (\cite{West}). 
We call it an {\it elliptic}, {\it hyperbolic} and {\it parabolic crosscap}, 
when $\alpha>0$,  $\alpha<0$,  $\alpha=0$, respectively.   
Obviously,  
when projecting the crosscap 
along its image tangent line $\ell=d\iota(TM_{x_0})$, 
the germ of $\Pi_\ell \circ \iota$ at $x_0$ is of corank $2$. 

\begin{thm} {\rm ({\bf Parabolic curve} \cite[Chap.5]{West})}  \label{west_thm}
The parabolic curve does not approach to any hyperbolic crosscap, 
while there are two smooth branches of parabolic curve 
approaching to any elliptic crosscap. 
For generic $\iota: M \to \R^3$ in the space of all maps having crosscaps, 
the singular germ $\Pi_\ell \circ \iota$ of corank $2$ 
is $\A$-equivalent to 
either the {\it sharksfin} 
or the {\it deltoid} in Table \ref{list2} 
(then we call it a generic elliptic or hyperbolic crosscap).  
\end{thm}

On the other hand, 
it seems that 
the flecnodal curve on the singular surface with crosscap 
had not been taken attention. 
In our previous paper, we showed that

\begin{thm}{\rm ({\bf Flecnodal curve} \cite[Thm. 1.3, 1.4]{YKO}\label{thmYKO})} 
The flecnodal curve does not approach to 
any hyperbolic crosscap, 
while there are two smooth branches of flecnodal curve 
approaching to any elliptic crosscap. 
Assume that $\iota$ has an elliptic crosscap at $x_0 \in M$. 
Then,  in the source space, 
each branch of the flecnodal curve is  
tangent to a branch of the parabolic curve  
with odd contact order at $x_0$; 
In particular, 
both pairs of branches have $3$-point contact 
if and only if 
the singular projection $\Pi_\ell \circ \iota$ is of type sharksfin 
(i.e, it is a generic elliptic crosscap). 
\end{thm}

\begin{figure}
\begin{center}
\includegraphics[width=11cm, clip]{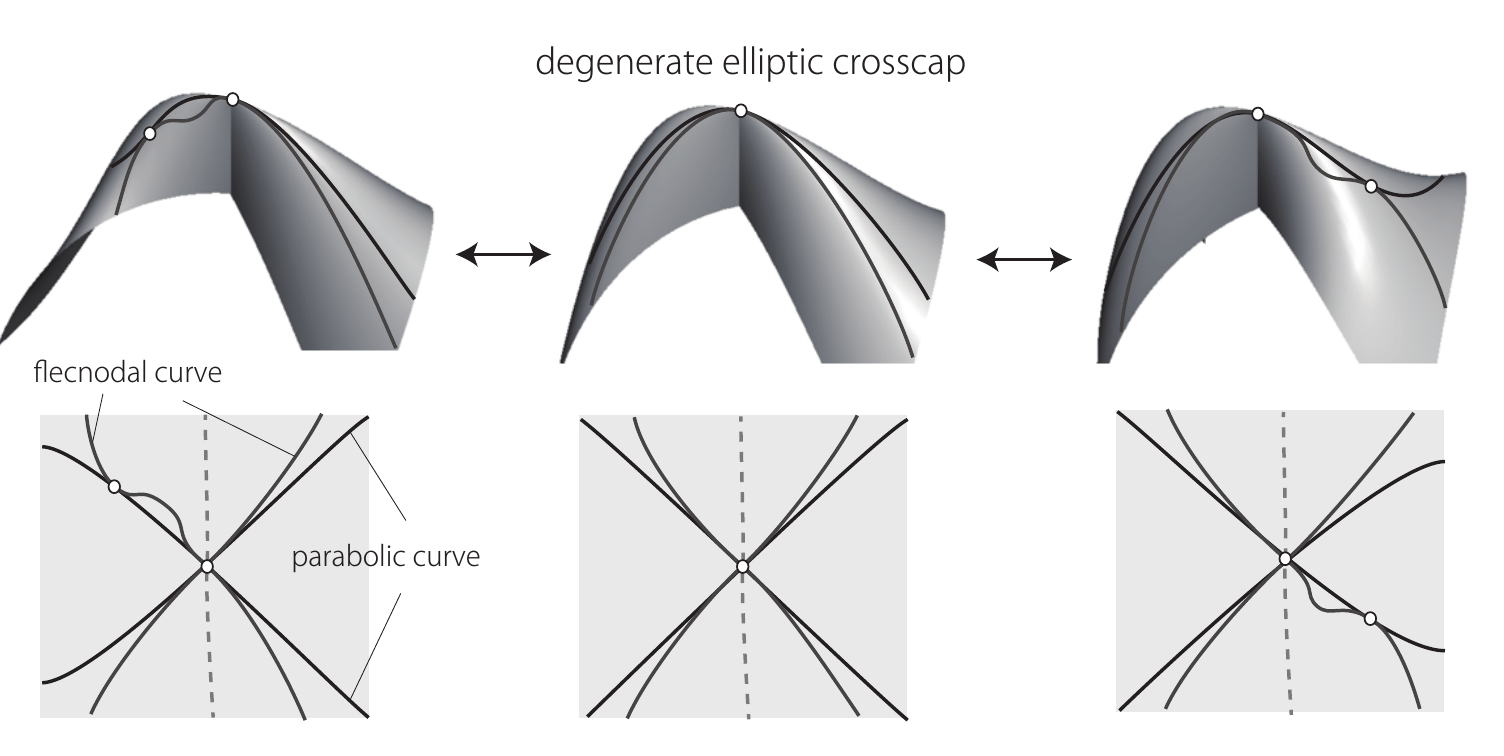}
\end{center}
\caption{\small   Bifurcation of degenerate elliptic crosscap \cite{YKO}. 
The upper pictures depict the bifurcation of the parabolic/flecnodal curves on 
the singular surface in $3$-space, 
while the lower pictures depict the corresponding loci in the source space of the parametrization 
(the dotted curve is the double point locus in the source). 
In the source, 
the degenerate crosscap (center)  has 
a pair of branches of parabolic/flecnodal curves with $5$-point contact. 
It is deformed into a generic elliptic crosscap (where two curves has $3$-point contact) 
and a godron ($2$-point contact). 
  }
\label{bif_elliptic}
\end{figure}

\begin{figure}
\begin{center}
\includegraphics[width=11cm, clip]{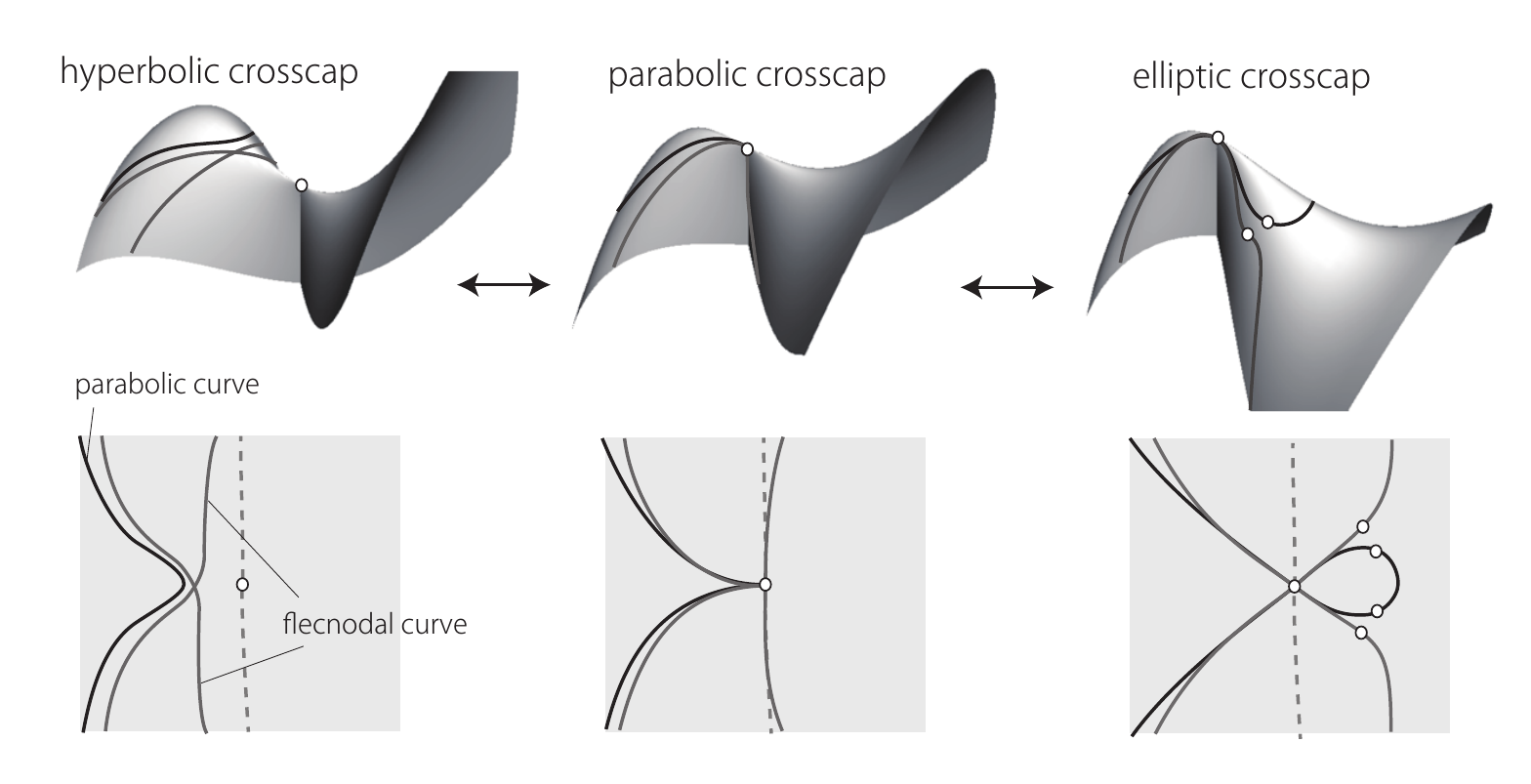}
\end{center}
\caption{\small   Bifurcation of parabolic crosscap. 
The bifurcation of the parabolic curve in source 
is given by a $1$-parameter family of sections of a crosscap 
(cf. Nu\~no-Ballesteros-Tari \cite{NT}). 
The flecnodal curve locally consists of two irreducible analytic branches. 
Deforming the parabolic crosscap, there appear 
a generic elliptic crosscap, two points of type goose, 
and two points of type butterfly, otherwise 
a hyperbolic crosscap and a double point of the flecnodal curve. 
The corresponding bifurcation diagram 
in the parameter space of view-lines $\ell$ 
is the same as  Fig.\ref{section_I23}.   }
\label{bif_parabolic}
\end{figure}

Now let us discuss a singular version of Rieger's observation 
on the projection of a moving surface \cite{Rieger2};   
Parallel projections of one-parameter families of crosscaps 
should be related to 
corank two map-germs of $\A_e$-codimension $3$. 
We obtain the following theorem:

\begin{thm}\label{thm} 
{\rm ({\bf Projection of non-generic crosscap})} 
For a generic one-parameter family 
$M\times I \to \R^3$, $(x,t) \mapsto \iota_t(x)$ 
of smooth maps having crosscaps (for all $t$), the germ of parallel projection 
$\Pi_\ell \circ \iota_0: M \to \R^3$ admits 
 the odd-shaped sharksfin and the type of $I_{2,3}$. 
The bifurcation of parabolic/flecnodal curves on the singular surface 
with respect to the parameter $t$ 
are described in Fig.\ref{bif_elliptic} and Fig.\ref{bif_parabolic}. 
\end{thm}

\proof 
In \cite{YKO} 
the assertion has been proved 
for the case of odd-shaped sharksfin, i.e., 
projecting the least non-generic elliptic crosscaps 
in the sense mentioned above.  
Here we deal with the case of $I_{2,3}$, 
i.e., projecting parabolic crosscaps. 

\begin{figure}
 \begin{center}
\includegraphics[width=4.5cm, clip]{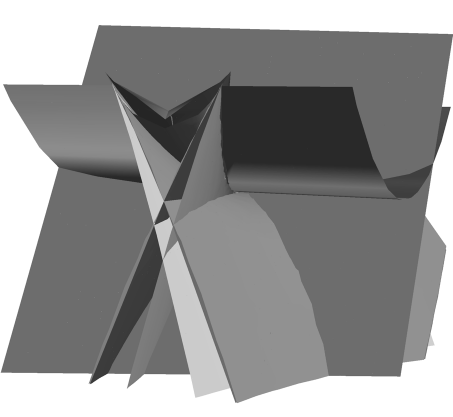}
\end{center}
\caption{\small   A generic plane section of $\B_G$ though the origin.   }
\label{section_I23_0}
\end{figure}

Suppose that we are given a smooth map 
$\iota: M \to \R^3$ with a parabolic crosscap at $x_0=(0,0)$.  
A generic $1$-parameter deformation 
$\iota_t: M \to \R^3$ of $\iota_0=\iota$ ($t$ sufficiently small) 
may be written as 
\begin{equation}\label{deform_crosscap}
\iota_t(x,y)=\left(y, xy+g(x,y,t), x^2+ t y^2+ \phi(x,y,t)\right),  
\end{equation}
where for each $t$ fixed, 
$j^2g(0)=0$, $j^2\phi(0)=0$ and $g(x,y,0)$ does not depend on $y$, 
via source coordinate changes 
and an affine transformation of target $\R^3$ depending on  $t$
(cf. \cite[Prop. 4.1]{Oliver}, \cite{NT}). 
The image tangent line at the parabolic crosscap 
is generated by $(1,0,0)$, and 
we put $U=\R^2$ of lines $\ell$ generated by $(1,v,w)$. 
The parallel projection of 
the singular surface $\iota_t(M)$ along $\ell$ defines  
a family of smooth maps $\varphi: M \times U \times \R \to \R^2$ by 
$$
\varphi(x, y, v, w, t)
=\left(xy-vy+g(x,y,t), x^2+t y^2-wy + \phi(x,y,t)\right). 
$$
For general choices of $g, \phi$ such as 
$\phi=c_{03}y^3+\cdots\; (c_{03}\not=0)$ etc, 
the plane-to-plane germ 
$\varphi(x,y,0,0,0)$ is $\K$-equivalent to $I_{2,3}$ 
and $\varphi(x,y,v,w,t)$ is a stable unfolding of this singularity type. 
Hence, generically, 
$\varphi$ can be regarded as a topologically versal unfolding 
of some germ belonging to the $\A$-moduli stratum of type $I_{2,3}$ (cf. Gaffney-Mond \cite{GM}). 
Thus the bifurcation diagram of $\varphi$ 
is obtained by a slightly modified $\B_G$ in Fig.\ref{I23} -- 
in particular, $\varphi(x,y,v,w,0)$ corresponds to a generic $2$-dimensional section 
of $\B_G$ through the origin, like as  in Fig.\ref{section_I23_0},  
and is homeomorphic to the middle of Fig.\ref{section_I23}. 
As $t$ varies from $0$, 
the section bifurcates in a  similar way as depicted 
in Fig.\ref{section_I23} (right and left). 

We have just seen the bifurcation of singular view-direction parameters $(u, v)$ 
with respect to $t$, and  now let us turn to see bifurcations of 
the parabolic curve and the flecnodal curve on the singular surface. 
Since we are interested in the topological type, 
we take at first a typical one of the form  (\ref{deform_crosscap}) 
$$
\iota_t(x,y)=\left(y, xy, x^2+ t y^2+y^3\right) 
$$
i.e., $g=0$, $\phi=y^3$. 
The  parallel projection $\varphi$ described above is equivalent to 
$G$ in (\ref{versal_I23}) via 
the coordinate changes $(\bar{x}, \bar{y})=(x-v, y)$ of the source 
and $(\bar{X}, \bar{Y})=(X, Y-v^2)$ of the target 
together with  $(a,b,c)=(2v, -w, t)$ of parameters. 
As shown in the proof of Theorem \ref{thm1}, 
the beaks-lips curve is given by 
$\bar{x}^2=\bar{y}^2(c+3\bar{y})$ 
and $a=-4\bar{x}$, 
and the swallowtail curve is given by 
$\bar{x}^2=\bar{y}^2(c+4\bar{y})$ 
and $a=\pm \bar{y}(4c+15\bar{y})(c+4\bar{y})^{-1/2}$. 
Hence substituting $(\bar{x}, \bar{y})=(x-a/2, y)$ to these equations, 
parabolic and flecnodal curves are obtained, respectively, such as 
$$\textstyle 
x^2=y^2(t+3y), \quad x^2(t+4y)=y^2(t+\frac{7}{2}y)^2.$$
When $t=0$, 
the parabolic curve has an ordinary cusp at the origin, 
while the flecnodal curve consists of two irreducible components, 
an ordinary cusp and a line $y=0$. 
The line is also the double point curve in this special from, 
but it is not the case in general: it can be seen that 
the component is tangent to the double-point curve at the origin. 
In fact, for general choices of $g, \phi$, 
above two equations are modified by adding 
functions of order $\ge 3$ (with respect to $x, y$), hence 
the local pictures around the origin do not change topologically, 
as depicted in Fig.\ref{bif_parabolic}.  \qed
 
 \begin{rem}{\rm 
 ({\bf Goose and butterfly}) 
 In Fig.\ref{bif_parabolic}, for $t>0$, 
the parabolic curve (black) 
has two parts corresponding to 
the beaks and the lips types of projection along 
the asymptotic line, which are separated by 
two points of goose type; 
the flecnodal curve (gray) has two points of 
type butterfly at which 
an asymptotic line has $5$-point contact with the surface. 
For $t<0$, 
the parabolic curve consists of a smooth component 
and an isolated point which is a hyperbolic crosscap, 
and 
the flecnodal curve consists of two smooth components 
meeting each other at a point transversely 
where both of two asymptotic lines have $5$-point contact 
with the surface. 
 }
 \end{rem}
 
\begin{rem}{\rm 
({\bf Foliation of asymptotic curves}) 
An {\it asymptotic curve} on a smooth surface in $\R^3$ is 
by definition the integral curve of the field of asymptotic lines, 
which is described  as the solution of 
a binary differential equation (Bruce-Tari \cite{BT}).  
All such curves form a pair of foliations on the hyperbolic domain 
and have cusp singularities at the parabolic points of the surface. 
The flecnodal curve is just the curve of inflection points of asymptotic curves. 
Near a crosscap point, 
the configuration of the asymptotic curves 
has been studied in Tari \cite{Tari}  (for elliptic and hyperbolic crosscaps) 
and Oliver \cite{Oliver} (for parabolic crosscaps) -- 
for instance, 
look at pictures of these foliations, 
Fig.1 and Fig.9 in  \cite{Oliver}, then 
we may experimentally trace the curve of inflection points of asymptotic curves, 
that would convince us of 
our pictures Fig.\ref{bif_elliptic} and Fig.\ref{bif_parabolic} above. 
We will discuss somewhere in detail 
about these two different approaches using BDE and parallel projection. 
}
\end{rem}

\section{Appendix: Planar caustics of parabolic umbilic type}

We briefly discuss an application to 
the planar parabolic umbilic caustics (Thom \cite{Thom}). 
First we explain a few definitions. 
Let us consider the $\mathcal{R}^+$-miniversal deformation of 
the $D_5$ singularity of function-germ
$$F(x,y, \mu_1, \mu_2, q_1, q_2):=x^2y+\frac{1}{4}y^4+\frac{q_2}{3}y^3+\frac{q_1}{2}y^2-2\mu_1x-\mu_2y$$
(usually one uses $x^2$ instead of the term of $y^3$). 
The catastrophe set $C$ in $\R^6=\R^2\times \R^4$ 
is defined by $\frac{\rd F}{\rd x}=\frac{\rd F}{\rd y}=0$, 
hence $C$ is parametrized by $x,y,q_1,q_2$ so that 
\begin{equation}\label{lag}
\mu_1= xy, \quad \mu_2=x^2+y^3+q_1y+q_2y^2.
\end{equation}
The corresponding {\it Lagrange map} $\Phi: C, 0 \to \R^4, 0$ 
is just the projection of $C$ 
to the parameter space, $(x,y,q_1, q_2) \mapsto (\mu_1, \mu_2, q_1, q_2)$. 

\begin{figure}
\begin{center}
\includegraphics[width=4.5cm]{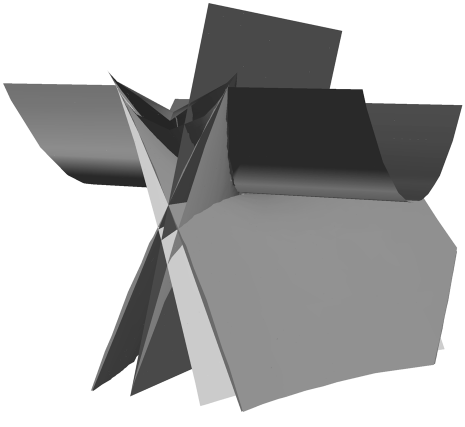}
\end{center}
\caption{\small   
Typical section of $\B_G$ with fixing the $c$-axis.  }
\label{plane_sections}
\end{figure}

\begin{figure}
\begin{center}
\includegraphics[width=11cm]{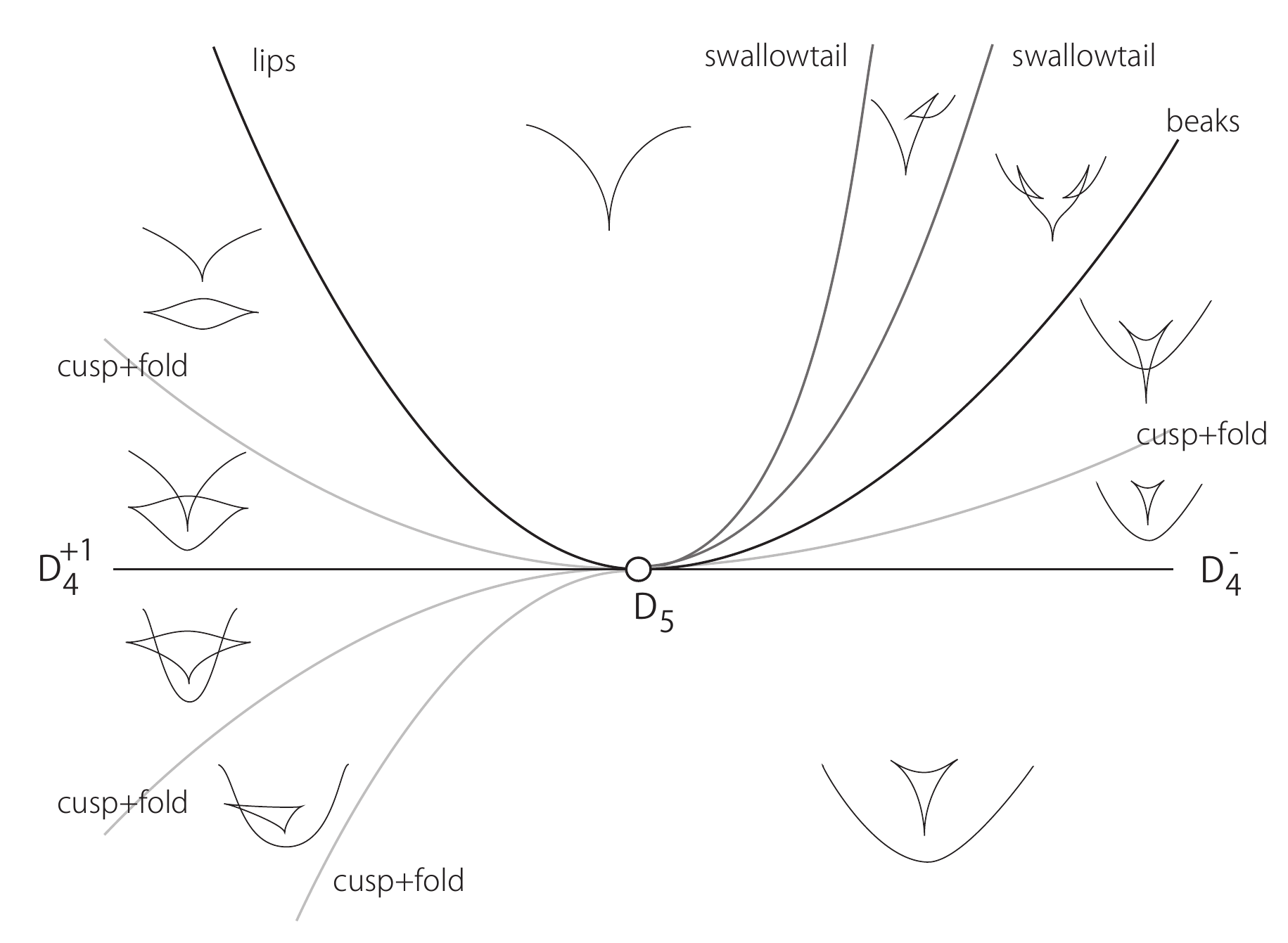}
\end{center}
\caption{\small   
Generic bifurcation of planar caustics of type $D_5$.  }
\label{D5}
\end{figure}

We are interested in 
generic $2$-dimensional sections of the big-caustics (the critical value set of $\Phi$) 
in $4$-space, see \cite[\S 4.1]{YKO} for a precise formulation. 
Let $\rho: \R^4,0 \to \R^2,0$  be a submersion 
so that $\Phi$ is transverse to 
the smooth level set $\rho^{-1}(0)$ at the origin. 
Note that if we put $S_t= \rho^{-1}(t)$ for each $t=(t_1, t_2) \in \R^2$ small enough, 
the condition says that $\Phi^{-1}(S_t)$ is smooth and 
the restrictions $\Phi: \Phi^{-1}(S_t) \to S_t$ 
form a $2$-parameter family of Lagrange maps. 
We now describe the family of Lagrange maps simply as a family of plane-to-plane maps 
(that corresponds to the so-called {\it caustics-equivalence}). 
Since $\Phi$ has corank $2$ and  $\rho\circ \Phi$ is submersive, 
we may assume that $\rho$ is of the form 
$$t_1=q_1-a_1\mu_1-a_2\mu_2+O(2), \quad 
t_2=q_2-b_1\mu_1-b_2\mu_2+O(2).$$
Solve $q_1, q_2$ in terms of $\mu_1, \mu_2, t_1, t_2$, and 
substitute them into (\ref{lag}), then  
by further coordinate changes of $(x,y)$ depending on $t_1, t_2$, we see that 
$\Phi$ is $\A$-equivalent to  
$\Xi: (x,y, t_1, t_2)\mapsto (\mu_1, \mu_2, t_1, t_2)$ with 
$$
\mu_1= xy, \quad 
\mu_2=x^2+y^3+a_1xy^2+ t_1y+t_2y^2 + \phi,$$  
where $j^3\phi(x,y,0,0)=0$ and $\frac{\rd \phi}{\rd t_i}(x,y,0,0)=0$. 
For general $\rho$,  we may think of $\Xi$ as an unfolding of the germ belonging to 
the $\A$-moduli of type $I_{2,3}$; 
Then the bifurcation diagram $\B_\Xi$ is obtained at least topologically as 
a small perturbation of the $bc$-plane section of $\B_G$ with fixing the $c$-axis, 
where $G$ is the normal form of (\ref{versal_I23}), see Fig.\ref{plane_sections}. 
Thus the topological type is unique as depicted in Fig.\ref{D5}. 

In $\R^4$ the big-caustics are only of type $A_\mu \; (\mu \le 5)$, $D_4$ and $D_5$. 
Each $\A$-type of corank one plane-to-plane germs 
can be realized by planar Lagrange map-germs, and hence 
generic multi-parameter bifurcations of planar caustics of type $A_\mu$ 
are the same as bifurcation diagrams of corank one germs. 
In \cite{YKO}, we have studied some topological types of 
 bifurcatoins of planar $D_4$-caustics, 
 and generic $2$-paramter bifurcation of type $D_5$ has been discussed above. 
Those lead a fairly natural extension of 
the {\it perestroikas} (=generic $1$-parameter bifurcations) 
of planar caustics due to Arnold-Zakalyukin \cite{Arnoldbooklet}. 

As an additional remark, 
each $\A$-type of corank two plane-to-plane germs 
can be realized by projecting singular Lagrange surface with open Whitney umbrella, 
see Bogaevski-Ishikawa \cite{BI}.

\begin{thebibliography}{99999}
\bibitem{Arnoldbooklet} V. I. Arnold, 
\textit{Catastrophe Theory}, 3rd edition, Springer (2004). 
\bibitem{Arnoldency} V. I. Arnold, V. V. Goryunov, O. V. Lyashko, V. A. Vasil'ev, 
\textit{Singularity Theory II, Classification and Applications}, 
Encyclopaedia of Mathematical Sciences Vol. 39, Dynamical System VIII (V. I. Arnold (ed.)), 
(translation from Russian version), Springer-Verlag  (1993). 
\bibitem{BI} I. A. Bogaevski and G. Ishikawa, 
Lagrange mappings of the first open Whitney umbrella, 
Pacific Jour. Math. 203 (1) (2002), 115--138. 
\bibitem{Bruce} J. W. Bruce, 
Projections and reflections of generic surfaces in $\R^3$, 
\textit{Math. Scand.} \textbf{54} (1984), 262--278. 
\bibitem{BT} J. W. Bruce and F. Tari, 
On binary differential equations, 
\textit{Nonlinearity} \textbf{8} (1995), 255--271.
\bibitem{BW} J. W. Bruce and J. West, 
Functions on a cross-cap, 
\textit{Math. Proc. Phil. Soc.} \textbf{123} (1998), 19--39. 
\bibitem{HF}  T. Fukui and M. Hasegawa, 
Singularities of parallel surfaces, 
\textit{Tohoku Math. J. } \textbf{64} (2012), 387-408.
\bibitem{GM} T.~Gaffney and D.~Mond, 
Weighted homogeneous maps from the plane to the plane, 
\textit{Math. Proc. Cambridge Phil. Soc.} \textbf{109} (1991), 451--470. 
\bibitem{GH}  C. G. Gibson and C. A. Hobbs, 
Singularity and Bifurcation for General Two Dimensional Planar Motions, 
\textit{New Zealand J. Math.} \textbf{25}  (1996) 141--163. 
\bibitem{Goryunov} V.V. Goryunov, 
Singularities of projections of complete intersections, 
\textit{J. Soviet Math.} \textbf{27} (1984), 2785--2811 
[\textit{Translated from Itogi Nauki i Tekhniki. Ser. Sovrem. Probl. Mat.} 
\textbf{22} (1983),167--206.]
\bibitem{HHUY} M. Hasegawa, A. Honda, K. Naokawa, M. Umehara and K. Yamada, 
Intrinsic invariants of cross caps, preprint, arXiv:1207.3853 (2012). 
\bibitem{Hawes} W. Hawes, Multi-dimensional Motions of the Plane and Space, 
\textit{Dissertation},  University of Liverpool (1994). 
\bibitem{Kabata} Y. Kabata, 
Recognition of plane-to-plane map-germs, 
preprint (2015), arXiv:1503.08544. 
\bibitem{Lander} L. Lander, The structure of the Thom-Boardman singularities of stable germs with type $\Sigma^{2,0}$, 
\textit{Proc. London Math. Soc.} \textbf{33} (1976), 113--137. 
\bibitem{NT} J. Nu\~no-Ballesteros and F. Tari, 
Surface in $\R^4$ and their projections to $3$-spaces, 
\textit{Proc. Royal Soc. Edinburgh Sect.} A \textbf{137} (2007), 1313--1328. 
\bibitem{OA}  T. Ohmoto and F. Aicardi, 
First order local  invariants of apparent contours, 
\textit{Topology}, \textbf{45} (2006) 27-45. 
\bibitem{Oliver} J. M. Oliver, 
On pairs of foliations of a parabolic cross-cap, 
\textit{Qual. Theory Dyn. Syst.}\textbf{10} (2011), 139--166. 
\bibitem{OT} R. Oset-Sinha and F. Tari, 
Projections of surfaces in $\R^4$ to $\R^3$ and geometry of their singular images, 
\textit{Rev. Math. Iberoam. European Math. Soc.} 
31 (2015), 33--50, DOI: 10.4171/RMI/825. 
\bibitem{Platonova} O. A. Platonova, 
Projections of smooth surfaces, 
\textit{J. Soviet Math.} \textbf{35} (1986), 2796--2808 
[\textit{Tr. Sem. I. G. Petrovskii} \textbf{10} (1984), 135--149 in Russian]. 
\bibitem{Rieger} J. H. Rieger, Families of maps from the plane to the plane, 
\textit{J. London Math. Soc.} (2) \textbf{36}  (1987), no. 2, 351-369.
\bibitem{Rieger2} J. H. Rieger, Versal topological stratification and the bifurcation geometry of map-germs of the plane, 
\textit{Math. Proc. Cambridge Philos.  Soc.}\textbf{107} (1990), 127-147.
\bibitem{Rieger3} J. H. Rieger, 
The geometry of view space of opaque objects bounded by smooth surfaces, 
\textit{Artificial Intelligence}, \textbf{44} (1990), 1-40.
\bibitem{RR}  J. H. Rieger and M. A. S. Ruas, Classification of $\A$-simple germs from $k^n$ to $k^2$, \textit{Compositio Math. } \textbf{79}  (1991), 99-108. 
\bibitem{Saji} K. Saji, 
Criteria for singularities of smooth maps from the plane into the plane and their applications, 
\textit{Hiroshima Math. Jour.} \textbf{40} (2010), 229--239. 
\bibitem{West} J. West, 
The Differential Geometry of the Cross-Cap, 
\textit{Dissertation}, University of Liverpool (1995). 
\bibitem{Tari} F. Tari, 
Pairs of geometric foliations on a crosscap, 
\textit{Tohoku Math. J.,} \textbf{59} (2007), 226--233.  
\bibitem{Thom} Thom,~R.,  
\textit{Structural Stability and Morphogenesis},
W. A. Benjam, (1972).
\bibitem{Yoshida} T. Yoshida, 
Bifurcation of plane-to-plane map-germs of corank $2$ and application to robotics (in Japanese), 
\textit{Master thesis}, Hokkaido University (2014).  
\bibitem{YKO} T. Yoshida, Y. Kabata and T. Ohmoto, 
Bifurcation of plane-to-plane map-germs of corank $2$, 
\textit{Quarterly Jour. Math.} (2014) doi:10.1093/qmath/hau013
\end{thebibliography}
\end{document}